\documentclass[a4paper,leqno,12pt]{amsart}
\usepackage[leqno]{amsmath}
\usepackage{amstext,amssymb,amsthm,enumitem}
\usepackage{graphicx}
\usepackage{tikz}
\usepackage[cp1250]{inputenc}
\usepackage{float}
\newtheorem*{theorem*}{Theorem}
\newtheorem{theorem}{Theorem}
\newtheorem{proposition}[theorem]{Proposition}
\theoremstyle{definition}
\newtheorem*{prof*}{Proof}
\newtheorem*{proof*}{Proof of Theorem 1}

\tikzset{dots/.append style={ultra thick, fill=none}}

\begin{document}
	
	\title[about weak and strong type inequalities for maximal operators]{On relations between weak and strong type inequalities for maximal operators on non-doubling metric measure spaces}
	
	\author{Dariusz Kosz}
	\address{ 
		\newline Faculty of Pure and Applied Mathematics
		\newline Wrocław University of Science and Technology 
		\newline Wyb. Wyspia\'nskiego 27 
		\newline 50-370 Wrocław, Poland
		\newline \textit{Dariusz.Kosz@pwr.edu.pl}	
	}

	\begin{abstract} In this article we characterize all possible cases that may occur in the relations between the sets of $p$ for which weak type $(p,p)$ and strong type $(p,p)$ inequalities for the Hardy--Littlewood maximal operators, both centered and non-centered, hold in the context of general metric measure spaces.
	
	
	\end{abstract}
	
	\maketitle

	\section{Introduction}
	Let $\mathbb{X} = (X, \rho, \mu)$ be a metric measure space with a metric $\rho$ and a Borel measure $\mu$ such that the measure of each ball is finite and strictly positive. Define the $\textit{Hardy--Littlewood}$ $\textit{maximal operators}$, centered $M^c$ and non-centered $M$, by
	\begin{displaymath}
	M^cf(x) = \sup_{r > 0} \frac{1}{\mu(B(x,r))} \int_{B(x,r)} |f| d\mu, \qquad x \in X,
	\end{displaymath}
	and
	\begin{displaymath}
	Mf(x) = \sup_{B \ni x} \frac{1}{\mu(B)} \int_B |f| d \mu , \qquad x \in X,
	\end{displaymath}
	respectively. Here $B$ refers to any open ball in $(X,\rho)$ and by $B(x,r)$ we denote the open ball centered at $x \in X$ with radius $r>0$.
	
	Recall that an operator $T$ is said to be of strong type $(p,p)$ for some $p \in [1, \infty]$ if $T$ is bounded on $L^p=L^p(\mathbb{X})$. Similarly, $T$ is of weak type $(p,p)$  if $T$ is bounded from $L^p$ to $L^{p,\infty}=L^{p,\infty}(\mathbb{X})$ (we use the convention $L^{\infty,\infty}=L^\infty$). Obviously, the operators $M^c$ and $M$ are of strong type $(\infty, \infty)$ in case of any metric measure space. Moreover, by using the Marcinkiewicz interpolation theorem, if $M^c$ (equivalently $M$) is of weak or strong type $(p_0,p_0)$ for some $p_0 \in [1,\infty)$, then it is of strong (and hence weak) type $(p,p)$ for every $p > p_0$. If the measure is doubling, that is $\mu(B(x, 2r)) \lesssim \mu(B(x, r))$ uniformly in $x \in X$ and $r > 0$, then both $M^c$ and $M$ are of weak type $(1,1)$. However, in general, the weak type $(1,1)$ inequalities may not occur. Furthermore, as we will see, it is even possible to construct a space for which the associated operators $M^c$ and $M$ are not of weak (and hence strong) type $(p,p)$ for every $p \in [1, \infty)$. 
	
	Finding examples of metric measure spaces with specific properties of associated maximal operators is usually a nontrivial task; see Aldaz \cite{Al}, for example. H.-Q. Li greatly contributed the program of searching spaces which are peculiar from the point of view of mapping properties of maximal operators. In this context, in \cite{Li1}, \cite{Li2} and \cite{Li3}, he considered a class of the cusp spaces. In \cite{Li1} H.-Q. Li showed, that for any fixed $1<p_0<\infty$ there exists a space for which the associated operator $M^c$ is of strong type $(p,p)$ if and only if $p > p_0$. Then, in \cite{Li2} examples of spaces were furnished for which $M$ is of strong type $(p,p)$ if and only if $p > p_0$. Moreover, for every $1 < \tau \leq 2$ there are examples of spaces for which $M^c$ is of weak type $(1,1)$, and $M$ is of strong type $(p,p)$ if and only if $p > \tau$. Finally, in \cite{Li3} H.-Q. Li showed that there are spaces with exponential volume growth for which $M^c$ is of weak type $(1,1)$, while $M$ is of strong type $(p,p)$ for every $p>1$. 
	
	The aim of this article is to complement and strengthen the results obtained by H.-Q. Li. For a fixed metric measure space $\mathbb{X}$ denote by $P_s^c$ and $P_s$ the sets consisting of such $p \in [1, \infty]$ for which the associated operators, $M^c$ and $M$ are of strong type $(p,p)$, respectively. Similarly, let $P_w^c$ and $P_w$ consist of such $p \in [1, \infty]$ for which $M^c$ and $M$ are of weak type $(p,p)$, respectively. Then\smallskip
	
	\begin{enumerate}[label=(\roman*)]
		\item each of the four sets is of the form $\{\infty\}$, $[p_0, \infty]$ or $(p_0,\infty]$, for some $p_0 \in [1, \infty)$;\smallskip
		
		\item we have the following inclusions
		\begin{displaymath}
		P_s \subset P_s^c, \quad P_w \subset P_w^c, \quad P_s^c \subset P_w^c \subset \overline{P_s^c}, \quad P_s \subset P_w \subset \overline{P_s},
		\end{displaymath}
		where $\overline{E}$ denotes the closure of $E$ in the usual topology of $\mathbb{R} \cup \{ \infty \}$.
	\end{enumerate}
	
	We will show that the conditions above are the only ones that the sets $P_s^c$, $P_s$, $P_w^c$ and $P_w$ must satisfy. Namely, we will prove the following.
	
	\begin{theorem}
		Let $P_s^c$, $P_s$, $P_w^c$ and $P_w$ be such that the conditions $(i)$ and $(ii)$ hold. Then there exists a (non-doubling) metric measure space for which the associated Hardy--Littlewood maximal operators, centered $M^c$ and non-centered $M$, satisfy\smallskip
		\begin{itemize}
			\item $M^c$ is of strong type $(p,p)$ if and only if $p \in P_s^c$,\smallskip
			\item $M$ is of strong type $(p,p)$ if and only if $p \in P_s$,\smallskip
			\item $M^c$ is of weak type $(p,p)$ if and only if $p \in P_w^c$,\smallskip
			\item $M$ is of weak type $(p,p)$ if and only if $p \in P_w$.\smallskip
		\end{itemize} 
	\end{theorem}
	The proof of Theorem 1 is postponed to Section 4.
	
	\section{First generation spaces}
	
	\smallskip We begin with a construction of some metric measure spaces called by us the $\textit{first} \linebreak \textit{generation spaces}$. The common property of these spaces is a similarity in the behavior of the associated operators $M^c$ and $M$, by what we mean the equalities $P_s^c = P_s$ and $P_w^c = P_w$. We begin with an overview of the first generation spaces and then we consider two subtypes separately in detail.
	
	Let $\tau = (\tau_n)_{n \in \mathbb{N}}$ be a fixed sequence of positive integers. Define 
	\begin{displaymath}
	X_{\tau} = \{x_n \colon n \in \mathbb{N}\} \cup \{x_{ni} \colon  i=1, \dots, \tau_n, n \in \mathbb{N}\},
	\end{displaymath}
	where all elements $x_n, x_{ni}$ are pairwise different (and located on the plane, say). We define the metric $\rho = \rho_\tau$ determining the distance between two different elements $x$ and $y$ by the formula
	\begin{displaymath}
	\rho(x,y) = \left\{ \begin{array}{rl}
	1 & \textrm{if } x_n \in \{x,y\} \subset S_n \textrm{ for some } n \in \mathbb{N},  \\
	2 & \textrm{in the other case.} \end{array} \right. 
	\end{displaymath}
	By $S_n$ we denote the branch $S_n = \{ x_n, x_{n1}, \dots , x_{n\tau_n} \}$ and by $S_n '$ the branch without the root, $S_n ' = S_n \setminus \{x_n\}$. Figure 1 shows a model of the space $(X_\tau, \rho)$. The solid line between two points indicates that the distance between them equals $1$. Otherwise the distance equals $2$.  
	
	\begin{figure}[H]
		\begin{tikzpicture}
		[scale=.8,auto=left,every node/.style={circle,fill,inner sep=2pt}]
		\node[label={[yshift=-1cm]$x_1$}] (n0) at (2.5,1) {};
		\node[label=$x_{11}$] (n1) at (1,4)  {};
		\node[label=$x_{12}$] (n2) at (2,4)  {};
		\node[label={[yshift=-0.06cm]$x_{1\tau_1}$}] (n3) at (4,4)  {};
		\node[dots] (n4) at (3,4)  {...};
		
		\node[label={[yshift=-1cm]$x_2$}] (m0) at (7.5,1) {};
		\node[label=$x_{21}$] (m1) at (6,4)  {};
		\node[label=$x_{22}$] (m2) at (7,4)  {};
		\node[label={[yshift=-0.06cm]$x_{2\tau_2}$}] (m3) at (9,4)  {};
		\node[dots] (m4) at (8,4)  {...};
		
		\node[dots] (o) at (11,1)  {...};
		
		\node[label={[yshift=-1cm]$x_n$}] (p0) at (14.5,1) {};
		\node[label=$x_{n1}$] (p1) at (13,4)  {};
		\node[label=$x_{n2}$] (p2) at (14,4)  {};
		\node[label={[yshift=-0.06cm]$x_{n\tau_n}$}] (p3) at (16,4)  {};
		\node[dots] (p4) at (15,4)  {...};
		
		\node[dots] (oo) at (18,1)  {...};
		
		\foreach \from/\to in {n0/n1, n0/n2, n0/n3, m0/m1, m0/m2, m0/m3, p0/p1, p0/p2, p0/p3}
		\draw (\from) -- (\to);
		\end{tikzpicture}
		\caption{}
	\end{figure}
	Note that we can explicitly describe any ball: for $n \in \mathbb{N}$,
	\begin{displaymath}
	B(x_n,r) = \left\{ \begin{array}{rl}
	\{x_n\} & \textrm{for } 0 < r \leq 1, \\
	S_n & \textrm{for } 1 < r \leq 2,  \\
	X_\tau & \textrm{for } 2 < r, \end{array} \right.
	\end{displaymath} 
	\noindent and for $i \in \{1, \dots, \tau_n\}$, $n \in \mathbb{N}$,
	\begin{displaymath}
	B(x_{ni},r) = \left\{ \begin{array}{rl}
	\{x_{ni}\} & \textrm{for } 0 < r \leq 1, \\
	\{x_n, x_{ni}\} & \textrm{for } 1 < r \leq 2,  \\
	X_\tau & \textrm{for } 2 < r. \end{array} \right.
	\end{displaymath} 
	
	We define the measure $\mu = \mu_{\tau, F}$ on $X_{\tau}$ by letting $\mu(\{x_n\}) = d_n$ and $\mu(\{x_{ni}\}) = d_nF(n,i)$, where $F > 0$ is a given function and $d = (d_n)_{n \in \mathbb{N}}$ is an appropriate sequence of strictly positive numbers with $d_1 = 1$ and $d_n$ chosen (uniquely!) in such a way that $\mu(S_n) = \mu(S_{n-1})/2$, $n \geq 2$. Note that this implies $\mu(X_\tau) < \infty$. Moreover, observe that $\mu$ is non-doubling. From now on we shall use the sign $|E|$ instead of $\mu(E)$ for $E \subset X_\tau$. It will be clear from the context when the symbol $|\cdot|$ refers to the measure and when it denotes the absolute value sign. 
	
	For a function $ f $ on $ X_\tau $ (which is in fact a 'sequence' of numbers) the Hardy--Littlewood maximal operators, centered $M^c$ and non-centered $M$, are given by
	\begin{displaymath}
	M^cf(x) = \sup_{r > 0} \frac{1}{|B(x,r)|} \sum_{y \in B(x,r)} |f(y)| \cdot |\{y\}|, \qquad x \in X_\tau,
	\end{displaymath}
	and
	\begin{displaymath}
	Mf(x) = \sup_{B \ni x} \frac{1}{|B|} \sum_{y \in B} |f(y)| \cdot |\{y\}|, \qquad x \in X_\tau,
	\end{displaymath}
	respectively. In this setting $M$ is of weak type $(p,p)$ for some $1 \leq p < \infty$ if $ \| M f \|_{p,\infty} \lesssim \| f \|_p$ uniformly in $f \in \ell^p(X_\tau, \mu)$, where $\| g \|_p = \big( \sum_{x \in X_\tau} |g(x)|^p |\{x\}|\big)^{1/p}$ and $ \| g \|_{p,\infty} = \sup_{\lambda > 0} \lambda |E_\lambda(g)|^{1/p}$, and $E_\lambda(g) = \{x \in X_\tau \colon |g(x)| > \lambda\}$. Similarly, $M$ is of strong type $(p,p)$ for some $1 \leq p \leq \infty$ if $ \| M f \|_p \lesssim \| f \|_p$ uniformly in $f \in \ell^p(X_\tau, \mu)$, where $\| g \|_\infty = \sup_{x \in X_\tau} |g(x)|$. Here the notation $A \lesssim B$ is used to indicate that $A \leq CB$ with a positive constant $C$ independent of significant quantities. Moreover, for given a function $f \geq 0$ and a set $E \subset X_\tau$ we denote the average value of $f$ on $E$ by 
	\begin{displaymath}
		A_E(f) = \frac{1}{|E|}\sum_{x \in E} f(x) |\{x\}|. 
	\end{displaymath}
	Analogous definitions and comments apply to $M^c$ instead of $M$ and then to both $M$ and $M^c$ in the context of the space $(Y_\tau, \mu)$ in Section 3.
	
	We are ready to describe two subtypes of the first generation spaces.   
	\subsection{} We first construct and investigate first generation spaces for which the equalities $P_s^c = P_s$ and $P_w^c = P_w$ hold and, in addition, there is no significant difference between the incidence of the weak and strong type inequalities, by what we mean that $P_s^c=P_w^c $ and $ P_s = P_w$. Of course, combining all these equalities, we obtain that for such spaces all four sets take the same form. In the first step, for any fixed $p_0 \in [1, \infty] $ we construct a space denoted by $\hat{\mathbb{X}}_{p_0}$ for which $P_s^c = P_s = P_w^c = P_w = [p_0, \infty]$ (by $[\infty, \infty]$ we mean $\{\infty\} $). Then, after slight modifications, for any fixed $p_0 \in [1, \infty)$ we get a space $\hat{\mathbb{X}}_{p_0}'$ for which $ P_s^c = P_s = P_w^c = P_w = (p_0, \infty]$.
	
	Fix $p_0 \in [1, \infty]$ and let $\hat{\mathbb{X}}_{p_0} = (X_\tau, \rho, \mu)$ be the first generation space with $\tau_n = \left \lfloor \frac{(n+1)^{p_0}}{n} \right \rfloor$ in the case $p_0 \in [1, \infty)$, or $\tau_n = 2^n$ in the case $p_0 = \infty$, and $F(n,i) = n$, $i=1, \dots \tau_n$, $n \in \mathbb{N}$. The key point for considerations that follow is that we have: for $p_0 \neq 1$,
	\begin{displaymath}
	\lim_{n \rightarrow \infty} \frac{n \tau_n}{(n+1)^p} = \infty, \qquad 1 \leq p < p_0, 
	\end{displaymath} 
	and for $p_0 \neq \infty$,
	\begin{displaymath}
	\frac{n \tau_n}{(n+1)^{p_0}} \leq 1, \qquad n \in \mathbb{N}.
	\end{displaymath}
	
	\begin{proposition}
		Fix $p_0 \in [1, \infty]$ and let $\hat{\mathbb{X}}_{p_0}$ be the metric measure space defined above. Then the associated maximal operators, centered $M^c$ and non-centered $M$, are not of weak type $(p,p)$ for $1 \leq p < p_0$, but are of strong type $(p,p)$ for $p \geq p_0$.
	\end{proposition}
	\begin{prof*}
		It suffices to prove that $M^c$ fails to be of weak type $(p,p)$ for $1 \leq p < p_0$ and $M$ is of strong type $(p_0,p_0)$. First we show that $M^c$ is not of weak type $(p,p)$ for $1 \leq p < p_0$. Consider $p_0 \in (1, \infty]$ and fix $p \in [1, p_0)$. Let $f_n = \delta_{x_n}$, $n \geq 1$. Then $\|f_n\|_p^p = d_n$ and $M^cf_n(x_{ni}) \geq \frac{1}{n+1}$, $i = 1, \dots, \tau_n$. This implies that $|E_{1/(2(n+1))}(M^cf_n)| \geq n \tau_n d_n$ and hence
		\begin{displaymath}
		\limsup_{n \rightarrow \infty} \frac{\|M^cf_n\|_{p,\infty}^p}{\|f_n\|_p^p} \geq \lim_{n \rightarrow \infty} \frac{n \tau_n d_n}{(2(n+1))^p d_n} = \infty.
		\end{displaymath} 
		
		In the next step we show that $M$ is of strong type $(p_0,p_0)$. Consider $p_0 \in [1, \infty)$, since the case $p_0 = \infty$ is trivial. Let $f \in \ell^{p_0}(\hat{\mathbb{X}}_{p_0})$. Without any loss of generality we assume that $f \geq 0$. Denote $\mathcal{D} = \{\{x_n, x_{ni}\} \colon n \in \mathbb{N}, i=1, \dots, \tau_n\}$. We use the estimate
		\begin{displaymath}
		\|Mf\|_{p_0}^{p_0} \leq \sum_{B \subset X_\tau} \sum_{x \in B} A_B(f)^{p_0} |\{x\}| = \sum_{B \subset X_\tau} A_B(f)^{p_0} |B|. 
		\end{displaymath}
		Note that each $x \in X_\tau$ belongs to at most three different balls which are not elements of $\mathcal{D}$. Combining this with Hölder's inequality, we obtain
		\begin{displaymath}
		\sum_{B \notin \mathcal{D}}A_B(f)^{p_0} |B| \leq \sum_{B \notin \mathcal{D}} \sum_{x \in B} f(x)^{p_0} |\{x\}| \leq 3 \|f\|_{p_0}^{p_0}. 
		\end{displaymath}
		Therefore
		\begin{equation} \label{2.1}
		\|Mf\|_{p_0}^{p_0} \leq 3 \|f\|_{p_0}^{p_0} + \sum_{n \in \mathbb{N}} \sum_{i=1}^{\tau_n} \Big(\frac{f(x_n)+ n f(x_{ni})}{n+1}\Big)^{p_0} |\{x_n, x_{ni}\}|.
		\end{equation}
		Finally, we use the inequalities $\big(f(x_n)+nf(x_{ni})\big)^{p_0} \leq \big(2f(x_n)\big)^{p_0} + \big(2nf(x_{ni})\big)^{p_0}$ and $|\{x_n, x_{ni}\}| \leq 2 |\{x_{ni}\}| = 2n |\{x_n\}|$ to estimate the double sum in (\ref{2.1}) by
		\begin{displaymath}
		2^{p_0+1} \Big( \sum_{n \in \mathbb{N}} \frac{n \tau_n }{(n+1)^{p_0}}f(x_n)^{p_0} |\{x_n\}| + \sum_{n \in \mathbb{N}} \sum_{i=1}^{\tau_n} \big(\frac{nf(x_{ni})}{n+1}\big)^{p_0} |\{x_{ni}\}|   \Big) \leq 2^{p_0+1} \|f\|_{p_0}^{p_0}. \raggedright \hfill \qed
		\end{displaymath}
		
	\end{prof*}
	
	A modification of arguments from the proof of Proposition 2 shows that for a fixed $p_0 \in [1, \infty)$, replacing the former $\tau_n$ by $\tau_n ' = \left \lfloor  (\log(n)+1) \frac{(n+1)^{p_0}}{n} \right \rfloor $, leads to the space $\hat{\mathbb{X}}_{p_0}'$ for which $P_s^c = P_s = P_w^c = P_w = (p_0, \infty]$. Moreover, it may be noted that only the properties $\lim_{n \rightarrow \infty}  \frac{n \tau_n '}{(n+1)^p} = \infty$, $1 \leq p \leq p_0$, and $\sup_{n \in \mathbb{N}} \frac{n \tau_n '}{(n+1)^p} < \infty$, $ p > p_0$, are essential.
	
	\subsection{} In contrast to the former case, for the spaces we now construct and study, the equalities $P_s^c = P_s$ and $P_w^c = P_w$ still hold, but there is a difference between the incidence of the weak and strong type inequalities. For any fixed $p_0 \in [1, \infty)$ we construct a space denoted by $\widetilde{\mathbb{X}}_{p_0}$ for which $P_s^c = P_s = (p_0, \infty]$ and $P_w^c = P_w = [p_0, \infty]$. We begin with the case $p_0 = 1$, which is discussed separately because it is relatively simple and may be helpful to outline the core idea behind the more difficult case $p_0 \in (1, \infty)$. 
	
	By $\widetilde{\mathbb{X}}_{1}$ we denote the first generation space $(X_\tau, \rho, \mu)$ with construction based on $\tau_n = n$ and $F(n,i) = 2^i$. Recall that $\mu$ is non-doubling.

	\begin{proposition}
		Let $\widetilde{\mathbb{X}}_{1}$ be the metric measure space defined above. Then the associated maximal operators, centered $M^c$ and non-centered $M$, are not of strong type $(1,1)$, but are of weak type $(1,1)$.
	\end{proposition}
	
	\begin{prof*}	
		First we note that $M^c$ fails to be of strong type $(1,1)$. Indeed, let $f_n = \delta_{x_n}$, $n \geq 1$. Then $\|f_n\|_1=d_n$ and for $i=1, \dots, n$ we have $M^cf_n(x_{ni}) \geq (1+2^i)^{-1} > 1/2^{i+1}$
		and hence $\|M^cf_n\|_1 \geq \sum_{i=1}^n 2^i d_n / 2^{i+1} = n \|f_n\|_1 / 2$.
		
		In the next step we show that $M$ is of weak type $(1,1)$. Let $f \in \ell^1(\widetilde{\mathbb{X}}_{1})$, $f \geq 0$, and consider $\lambda > 0$ such that $E_\lambda (Mf) \neq \emptyset$. If $\lambda < A_{X_\tau}(f)$, then $\lambda |E_\lambda(Mf)| / \|f\|_1 \leq 1$ follows. Therefore, from now on assume that $\lambda \geq A_{X_\tau}(f)$. With this assumption, if for some $x \in S_n$ we have $Mf(x) > \lambda$, then any ball $B$ containing $x$ and realizing $A_B(f) > \lambda$ must be a subset of $S_n$. Take any $n \in \mathbb{N}$ such that $E_\lambda(Mf) \cap S_n \neq \emptyset$. If $\lambda < A_{S_n}(f)$, then
		\begin{equation} \label{2.2.1}
		\frac{\lambda |E_\lambda(Mf) \cap S_n|}{\sum_{x \in S_n} f(x)|\{x\}|} \leq 1.
		\end{equation}
		Assume that $\lambda \geq A_{S_n}(f)$ and take $x \in E_\lambda(Mf) \cap S_n$. Now, any ball $B$ containing $x$ and realizing $A_B(f) > \lambda$ must be a proper subset of $S_n$. If $E_\lambda(Mf) \cap S_n' = \emptyset$, then $x = x_n$ so we obtain $f(x_n) > \lambda$ and hence (\ref{2.2.1}) again follows. In the opposite case, if $E_\lambda(Mf) \cap S_n' \neq \emptyset$, denote $j = \max\{i\in \{1, \dots, n\} \colon Mf(x_{ni})> \lambda\}$. Then $f(x_{nj}) > \lambda$ or $\frac{f(x_n)|\{x_n\}| + f(x_{nj}) |\{x_{nj}\}|}{|\{x_n\}| + |\{x_{nj}\}|} > \lambda$. Therefore, $f(x_n)|\{x_n\}| + f(x_{nj}) |\{x_{nj}\}| > \lambda |\{x_{nj}\}|$ and combining this with the estimate $|E_\lambda(Mf) \cap S_n| \leq 2|\{x_{nj}\}|$, we obtain
		\begin{displaymath}
		\frac{\lambda |E_\lambda(Mf) \cap S_n|}{\sum_{x \in S_n} f(x)|\{x\}|} \leq \frac{2\lambda|\{x_{nj}\}|}{f(x_n)|\{x_n\}| + f(x_{nj}) |\{x_{nj}|\}} \leq 2.
		\end{displaymath}
		
		Since $\frac{\lambda |E_\lambda(Mf) \cap S_n|}{\sum_{x \in S_n} f(x)|\{x\}|} \leq 2$ for any branch $S_n$ such that $E_\lambda(Mf) \cap S_n \neq \emptyset$, we have
		\begin{displaymath}
		\frac{\lambda |E_\lambda(Mf)|}{\|f\|_1} \leq 2,
		\end{displaymath}
		and, consequently, the weak type $(1,1)$ estimate $\|Mf\|_{1, \infty} \leq 2 \|f\|_{1}$ follows.
		$\raggedright \hfill \qed$ 
	\end{prof*}
	
	Now fix $p_0 \in (1, \infty)$ and consider $\widetilde{\mathbb{X}}_{p_0} = (X_\tau, \rho, \mu)$, with construction based on $\tau_n = \tau_n(p_0)$ and $F(n,i) = F_{p_0}(n,i)$, defined as follows. Let $c_n = \left \lfloor \frac{(n+1)^{p_0}}{n} \right \rfloor $ and 
	\begin{displaymath}
	e_n = \max\{k \in \mathbb{N} \colon 2^{k-1} \leq c_n \text{ and } 2^{1 - k - p_0} \geq \big(\frac{1}{1+n}\big)^{p_0}\}, \quad n \in \mathbb{N}.
	\end{displaymath}
	Note that $\lim_{n \rightarrow \infty} e_n = \infty$. Then, for $j \in \{1, \dots, e_n\}$, $n \in \mathbb{N}$, define $m_{nj}$ by the equality
	\begin{displaymath}
	2^{1-j} \Big( \frac{1}{1+m_{nj}} \Big)^{p_0}  = \Big( \frac{1}{1+n} \Big)^{p_0},
	\end{displaymath}
	and $s_{nj}$ by
	\begin{displaymath}
	s_{nj}= \min\{k \in \mathbb{N} \colon km_{nj} \geq 2^{1-j}nc_n\}.
	\end{displaymath}
	Observe that for $j \in \{1, \dots, e_n\}, n \in \mathbb{N}$,
	\begin{displaymath}
	1 \leq m_{nj} \leq n, \quad 2^{1-j}nc_n \leq s_{nj}m_{nj} \leq 2^{2-j}nc_n. 
	\end{displaymath}
	Finally, denote $\tau_n = \sum_{j=1}^{e_n} s_{nj}$, $n \in \mathbb{N}$, and $F(n,i) = m_{nj(n,i)}$, $i =1, \dots, \tau_n$, $n \in \mathbb{N}$, where
	\begin{displaymath}
	j(n,i) = \min\{k \in \{1, \dots, e_n\} \colon \sum_{j=1}^k s_{nj} \geq i\}.
	\end{displaymath}
	
	\begin{proposition}
		Let $\widetilde{\mathbb{X}}_{p_0}$ be the metric measure space defined above. Then the associated maximal operators, centered $M^c$ and non-centered $M$, are not of strong type $(p_0,p_0)$, but are of weak type $(p_0,p_0)$.
	\end{proposition}
	
	\begin{prof*}
		First we note that $M^c$ is not of strong type $(p_0,p_0)$. Indeed, let $f_n = \delta_{x_n}$, $n \geq 1$. Then $\|f_n\|_{p_0}^{p_0}=d_n$ and for $i=1, \dots, \tau_n$ we have $M^cf_n(x_{ni}) \geq (1+m_{nj(n,i)})^{-1}$ and hence
		\begin{align*}
		\|M^cf_n\|_{p_0}^{p_0} & \geq \sum_{j=1}^{e_n} \sum_{k=1}^{s_{nj}}\Big(\frac{1}{1+m_{nj}}\Big)^{p_0} d_n m_{nj} = d_n \sum_{j=1}^{e_n} \frac{s_{nj}m_{nj}}{(1+m_{nj})^{p_0}} \\ & \geq d_n \sum_{j=1}^{e_n} \frac{2^{1-j} n c_n}{(1+m_{nj})^{p_0}} =  d_n \sum_{j=1}^{e_n}\frac{nc_n}{(1+n)^{p_0}} = e_n \frac{nc_n}{(1+n)^{p_0}} \|f_n\|_{p_0}^{p_0}.
		\end{align*}
		Since $\lim_{n \rightarrow \infty} e_n = \infty$ and $\lim_{n \rightarrow \infty} \frac{n c_n}{(1+n)^{p_0}} = 1$, we are done.
		
		In the next step we show that $M$ is of weak type $(p_0,p_0)$. Let $f \in \ell^{p_0}(\widetilde{\mathbb{X}}_{p_0})$, $f \geq 0$, and consider $\lambda > 0$ such that $E_\lambda (Mf) \neq \emptyset$. If $\lambda < A_{X_\tau}(f)$, then using the inequality $\|f\|_1 \leq \|f\|_{p_0} |X_\tau|^{1/q_0}$, where $q_0$ is the exponent conjugate to $p_0$, we obtain $\lambda^{p_0} |E_\lambda(Mf)| / \|f\|_{p_0}^{p_0} <1$. Therefore, from now on assume that $\lambda \geq A_{X_\tau}(f)$. Take any $n \in \mathbb{N}$ such that $E_\lambda(Mf) \cap S_n \neq \emptyset$. If $\lambda < A_{S_n}(f)$, then using similar argument as above we have
		\begin{equation} \label{2.2.2}
		\frac{\lambda^{p_0} |E_\lambda(Mf) \cap S_n|}{\sum_{x \in S_n} f(x)^{p_0} |\{x\}|} \leq 1.
		\end{equation}
		Assume that $\lambda \geq A_{S_n}(f)$. If $E_\lambda(Mf) \cap S_n' = \emptyset$, then $f(x_n) > \lambda$ and hence (\ref{2.2.2}) again follows. In the opposite case, we have $|E_\lambda(Mf) \cap S_n| \leq 2 |E_\lambda(Mf) \cap S_n '|$. Assume that $f(x_n)< (1+m_{ne_n})\lambda / 2$. If $x \in E_\lambda(Mf) \cap S_n'$, then $f(x) \geq \lambda / 2$ and hence
		\begin{displaymath}
		\frac{\lambda^{p_0} |E_\lambda(Mf) \cap S_n|}{\sum_{x \in S_n} f(x)^{p_0} |\{x\}|} \leq \frac{2 \lambda^{p_0} |E_\lambda(Mf) \cap S_n'|}{\sum_{x \in S_n} f(x)^{p_0} |\{x\}|} \leq 2^{p_0+1}.
		\end{displaymath}
		Otherwise, if $f(x_n) \geq (1+m_{ne_n})\lambda / 2$, denote $r = \min\{j \in \{1, \dots, e_n\} \colon f(x_n) \geq (1+m_{nj})\lambda / 2 \}$. Let $S_n^{(r)} = \{x_{ni} \colon i \in \{1, \dots, \sum_{j=1}^{r-1}s_{nj}\}\}$. Note that the case $S_n^{(r)} = \emptyset$ is possible. Assume that $S_n^{(r)} \neq \emptyset$. If $x \in E_\lambda(Mf) \cap S_n^{(r)}$, then $f(x) > \lambda / 2$ and hence
		\begin{displaymath}
		\frac{\lambda^{p_0} |E_\lambda(Mf) \cap S_n^{(r)}|}{\sum_{x \in S_n^{(r)}} f(x)^{p_0} |\{x\}|} \leq 2^{p_0+1}.
		\end{displaymath}
		Moreover, we have
		\begin{align*}
		\frac{\lambda^{p_0} |E_\lambda(Mf) \cap \big(S_n \setminus S_n^{(r)}\big)|}{f(x_n)^{p_0} |\{x_n\}|} & \leq \Big(\frac{2}{1+m_{nr}}\Big)^{p_0} \frac{|S_n \setminus S_n^{(r)}|}{|\{x_n\}|} \leq \Big(\frac{2}{1+m_{nr}}\Big)^{p_0} \frac{2 |(S_n \setminus S_n^{(r)}) \cap S_n'|}{|\{x_n\}|} \\
		& \leq \Big(\frac{2}{1+m_{nr}}\Big)^{p_0} 2 \sum_{j=r}^{e_n}n c_n 2^{2-j} < 2^{p_0+4-r} n c_n \Big(\frac{1}{1+m_{nr}}\Big)^{p_0} \\ & = 2^{p_0+3} \frac{nc_n}{(1+n)^{p_0}} \leq 2^{p_0+3}.
		\end{align*}
		Therefore, regardless of the posibilities, $S_n^{(r)} = \emptyset$ or $S_n^{(r)} \neq \emptyset$, we obtain $\frac{\lambda^{p_0} |E_\lambda(Mf) \cap S_n|}{\sum_{x \in S_n} f(x)^{p_0} |\{x\}|} \leq 2^{p_0+3}$. 
		Since $\lambda^{p_0} |E_\lambda(Mf) \cap S_n| / \sum_{x \in S_n} f(x)^{p_0} |\{x\}| \leq 2^{p_0+3}$ for any branch $S_n$ such that $E_\lambda(Mf) \cap S_n \neq \emptyset$, we have $\lambda^{p_0} |E_\lambda(Mf)| / \|f\|_{p_0}^{p_0} \leq 2^{p_0+3}$ and, consequently, $\|Mf\|_{p_0, \infty}^{p_0} \leq 2^{p_0+3} \|f\|_{p_0}^{p_0}$.
		$\raggedright \hfill \qed$
	\end{prof*}
	\section{Second generation spaces}
	Now we construct and study metric measure spaces called by us the $\textit{second generation}$ $\textit{spaces}$. The common attribute of these spaces is a significant difference in the behavior of the associated operators $M^c$ and $M$, by what we mean that $M^c$ is of strong type $(1,1)$, which implies the basic property $P_s^c = P_w^c = [1, \infty]$, while $P_s$ (and possibly $P_w$) is a proper subset of $[1, \infty]$. Let $\tau = (\tau_n)_{n \in \mathbb{N}}$ be a fixed sequence of positive integers. Define
	\begin{displaymath}
	Y_\tau = \{y_n \colon n \in \mathbb{N}\} \cup \{y_{ni}, y_{ni}' \colon i = 1, \dots, \tau_n, n \in \mathbb{N} \},
	\end{displaymath}
	where all elements $y_n, y_{ni}, y_{ni}'$ are pairwise different. We define the metric $\rho = \rho_\tau$ determining the distance between two different elements $x$ and $y$ by the formula
	\begin{displaymath}
	\rho(x,y) = \left\{ \begin{array}{rl}
	1 & \textrm{if } \{x,y\} = T_{ni} \textrm{ or } y_n \in \{x,y\} \subset T_n \setminus T_n ' \textrm{ for some } n \in \mathbb{N}, \ i \in \{1, \dots, \tau_n\}, \\
	2 & \textrm{in the other case.} \end{array} \right. 
	\end{displaymath}
	By $T_n$ we denote the branch $T_n = \{ y_n, y_{n1}, \dots , y_{n\tau_n}, y_{n1}', \dots , y_{n\tau_n}' \}$. Additionally, we denote $T_n ' = \{y_{n1}', \dots , y_{n\tau_n}'\}$ and $T_{ni}=\{y_{ni}, y_{ni}'\}$. Figure 2 shows a model of the space $(Y_\tau, \rho)$.
	\begin{figure}[H]
		\begin{tikzpicture}
		[scale=.8,auto=left,every node/.style={circle,fill,inner sep = 2pt}]
		\node[label={[yshift=-1cm]$y_1$}] (n0) at (4,1) {};
		\node[label=$y_{11}$] (n1) at (2.5,3)  {};
		\node[label={[xshift=0.3cm]$y_{12}$}] (n2) at (3.5,3)  {};
		\node[label={[yshift=-0.06cm]$y_{1\tau_1}$}] (n3) at (5.5,3)  {};
		\node[dots] (n4) at (4.5,3)  {...};
		\node[label=$y_{11}'$] (n5) at (1,5)  {};
		\node[label=$y_{12}'$] (n6) at (3,5)  {};
		\node[label={[yshift=-0.06cm]$y_{1\tau_1}'$}] (n7) at (7,5)  {};
		\node[dots] (n8) at (5,5)  {...};
		
		\node[dots] (n8) at (9,1)  {...};
		
		\node[label={[yshift=-1cm]$y_n$}] (m0) at (14,1) {};
		\node[label=$y_{n1}$] (m1) at (12.5,3)  {};
		\node[label={[xshift=0.3cm]$y_{n2}$}] (m2) at (13.5,3)  {};
		\node[label={[yshift=-0.06cm]$y_{n\tau_n}$}] (m3) at (15.5,3)  {};
		\node[dots] (m4) at (14.5,3)  {...};
		\node[label=$y_{n1}'$] (m5) at (11,5)  {};
		\node[label=$y_{n2}'$] (m6) at (13,5)  {};
		\node[label={[yshift=-0.06cm]$y_{n\tau_n}'$}] (m7) at (17,5)  {};
		\node[dots] (m8) at (15,5)  {...};
		
		\node[dots] (n8) at (19,1)  {...};
		
		\foreach \from/\to in {n0/n1, n0/n2, n0/n3, n1/n5, n2/n6, n3/n7, m0/m1, m0/m2, m0/m3, m1/m5, m2/m6, m3/m7}
		\draw (\from) -- (\to);
		\end{tikzpicture}
		\caption{}
	\end{figure}
	Note that we can explicitly describe any ball: for $n \in \mathbb{N}$,
	\begin{displaymath}
	B(y_n,r) = \left\{ \begin{array}{rl}
	\{y_n\} & \textrm{for } 0 < r \leq 1, \\
	T_n \setminus T_n ' & \textrm{for } 1 < r \leq 2,  \\
	Y_\tau & \textrm{for } 2 < r, \end{array} \right.
	\end{displaymath} 
	\noindent and for $i \in \{1, \dots, \tau_n\}$, $n \in \mathbb{N}$,
	\begin{displaymath}
	B(y_{ni},r) = \left\{ \begin{array}{rl}
	\{y_{ni}\} & \textrm{for } 0 < r \leq 1, \\
	\{y_n\} \cup T_{ni} & \textrm{for } 1 < r \leq 2,  \\
	Y_\tau & \textrm{for } 2 < r, \end{array} \right. \qquad
	\end{displaymath} 
	and
	\begin{displaymath}
	B(y_{ni}',r) = \left\{ \begin{array}{rl}
	\{y_{ni}'\} & \textrm{for } 0 < r \leq 1, \\
	T_{ni} & \textrm{for } 1 < r \leq 2,  \\
	Y_\tau & \textrm{for } 2 < r. \end{array} \right. \qquad
	\end{displaymath}
	
	We define the measure $\mu = \mu_{\tau, F}$ by letting $\mu(\{y_n\}) = d_n$, $\mu(\{y_{ni}\}) = \frac{d_n}{\tau_n}$ and $\mu(\{y_{ni}'\}) = d_n
	F(n,i)$, where $F > 0$ is a given function and $d = (d_n)_{n \in \mathbb{N}}$ is an appropriate sequence of strictly positive numbers with $d_1 = 1$ and $d_n$ chosen (uniquely!) in such a way that $|T_n| = |T_{n-1}|/2$, $n \geq 2$. Note that this implies $|Y_\tau| < \infty$ and observe that $\mu$ is non-doubling.
	
	We are ready to describe two subtypes of the second generation spaces.
	
	\subsection{} We first construct spaces for which apart from the basic property $P_s^c = P_w^c = [1, \infty]$ we also have $P_s = P_w$. In the first step, for any fixed $p_0 \in (1, \infty]$ we construct a space denoted by $\hat{\mathbb{Y}}_{p_0}$ for which $P_s = P_w = [p_0, \infty]$. Then, after a slight modification, for any fixed $p_0 \in [1, \infty)$ we get a space $\hat{\mathbb{Y}}_{p_0}'$ for which $P_s = P_w = (p_0, \infty]$.
	
	Fix $p_0 \in (1, \infty]$ and let $\hat{\mathbb{Y}}_{p_0}$ be the second generation space with $ \tau_n = \left \lfloor \frac{(n+1)^{p_0}}{n} \right \rfloor$ in the case $p_0 \in (1, \infty)$, or $\tau_n = 2^n$ when $p_0 = \infty$, and $F(n,i) = n$, $i = 1, \dots, \tau_n$, $n \in \mathbb{N}$. 
	
	\begin{proposition}
		Let $\hat{\mathbb{Y}}_{p_0}$ be the metric measure space defined above. Then the associated centered maximal operator $M^c$ is of strong type $(1,1)$, while the non-centered $M$ is not of weak type $(p,p)$ for $1 \leq p < p_0$, but is of strong type $(p,p)$ for $p \geq p_0$.
	\end{proposition}
	
	\begin{prof*}
		First we show that $M^c$ is of strong type $(1,1)$. Let $f \in \ell^1(\hat{\mathbb{Y}}_{p_0})$, $f \geq 0$. Denote $\mathcal{G} = \{\{y_n\} \cup T_{ni} \colon n \in \mathbb{N}, i = 1, \dots, \tau_n\}$ and $\mathcal{B}_y = \{B(y, \frac{1}{2}),B(y, \frac{3}{2}), B(y, \frac{5}{2})\}$, $y \in Y_\tau$. We use the estimate
		\begin{displaymath}
		\|M^cf\|_1 \leq \sum_{y \in Y_\tau} \sum_{B \in \mathcal{B}_y} A_{B}(f)  |\{y\}|.
		\end{displaymath}
		Note that each $y \in Y_\tau$ belongs to at most four different balls which are not elements of $\mathcal{G}$. Thus we obtain
		\begin{displaymath}
		\sum_{y \in Y_\tau} \sum_{B \in \mathcal{B}_y \setminus \mathcal{G}} A_{B}(f)  |\{y\}| \leq \sum_{B \notin \mathcal{G}} \sum_{y \in B} f(y)|\{y\}| \leq 4 \|f\|_1.
		\end{displaymath}
		Therefore
		\begin{displaymath}
		\|M^cf\|_1 \leq 4 \|f\|_1 + \sum_{n \in \mathbb{N}} \sum_{i=1}^{\tau_n} A_{B(y_{ni}, \frac{3}{2})}(f)  |\{y_{ni}\}|.
		\end{displaymath}
		It suffices to see that the last term of the above expression is estimated by
		\begin{displaymath}
		\sum_{n \in \mathbb{N}} \tau_n f(y_n) |\{y_{n1}\}| + \sum_{n \in \mathbb{N}} \sum_{i=1}^{\tau_n} \Big( f(y_{ni}) |\{y_{ni}\}| + f(y_{ni}') |\{y_{ni}'\}| \Big) = \|f\|_1.
		\end{displaymath}
		
		In the next step we show that $M$ is not of weak type $(p,p)$ for $1 \leq p < p_0$. Indeed, fix $p < p_0$ and let $f_n = \delta_{y_n}$, $n \geq 1$. Then $\|f_n\|_p^p = d_n$ and $Mf_n(y_{ni}') \geq \frac{1}{n+1+(1/\tau_n)}\geq \frac{1}{n+2}$, $i = 1, \dots, \tau_n$. This implies that $|E_{1/(2(n+2))}(M f_n)| \geq n \tau_n d_n$ and hence we obtain
		\begin{displaymath}
		\limsup_{n \rightarrow \infty} \frac{\|Mf_n\|_{p,\infty}^p}{\|f_n\|_p^p} \geq \lim_{n \rightarrow \infty} \frac{n \tau_n d_n}{(2(n+2))^p d_n} = \infty.
		\end{displaymath} 
		
		To complete the proof, it suffices to show that $M$ is of strong type $(p_0,p_0)$ in the case $p_0 \in (1, \infty)$. Let $f \in \ell^{p_0}(\hat{\mathbb{Y}}_{p_0})$, $f \geq 0$. We use the estimate
		\begin{displaymath}
		\|Mf\|_{p_0}^{p_0} \leq \sum_{B \subset Y_\tau} \sum_{y \in B} A_B(f)^{p_0} |\{y\}| = \sum_{B \subset Y_\tau} A_B(f)^{p_0} |B|. 
		\end{displaymath}
		Once again note that each $y \in Y_\tau$ belongs to at most four different balls which are not elements of $\mathcal{G}$. Combining this with Hölder's inequality, we obtain
		\begin{displaymath}
		\sum_{B \notin \mathcal{G}} A_B(f)^{p_0} |B| \leq \sum_{B \notin \mathcal{G}} \sum_{y \in B} f(y)^{p_0} |\{y\}| \leq 4 \|f\|_{p_0}^{p_0}. 
		\end{displaymath}
		Therefore
		\begin{equation} \label{3.1}
		\|Mf\|_{p_0}^{p_0} \leq 4 \|f\|_{p_0}^{p_0} + \sum_{n \in \mathbb{N}} \sum_{i=1}^{\tau_n} \Big(\frac{f(y_n)+ 1/{\tau_n}f(y_{ni}) + n f(y_{ni}')}{1 + 1/{\tau_n}+n}\Big)^{p_0} |\{y_n, y_{ni}, y_{ni}'\}|.
		\end{equation}
		Finally, we use the inequalities 
		\begin{displaymath}
		\big(f(y_n)+ 1/{\tau_n}f(y_{ni}) + n f(y_{ni}')\big)^{p_0} \leq \big(3f(y_n)\big)^{p_0} + \big(3f(y_{ni})/{\tau_n}\big)^{p_0} + \big(3nf(y_{ni}')\big)^{p_0},
		\end{displaymath}
		and $|\{y_n, y_{ni}, y_{ni}'\}| \leq 3 |\{y_{ni}'\}| = 3n |\{y_n\}|$ to estimate the double sum in (\ref{3.1}) by
		\begin{displaymath}
		3^{p_0+1} \Big( \sum_{n \in \mathbb{N}} \frac{n \tau_n f(y_n)^{p_0} }{(n+1)^{p_0}} |\{y_n\}| + \sum_{n \in \mathbb{N}} \sum_{i=1}^{\tau_n} \frac{\big(f(y_{ni})/{\tau_n}\big)^{p_0} + (nf(y_{ni}')\big)^{p_0}}{(1+1/{\tau_n}+n)^{p_0}} |\{y_{ni}'\}| \Big) \leq 3^{p_0+1} \|f\|_{p_0}^{p_0}.
		\end{displaymath}
		$\raggedright \hfill \qed$
	\end{prof*}
	
	Note that in the same way as it was done at the end of Section 2.1, replacing the former $\tau_n$ by $\tau_n ' = \left \lfloor  (\log(n)+1) \frac{(n+1)^{p_0}}{n} \right \rfloor $, $p_0 \in [1, \infty)$, results in obtaining the space $\hat{\mathbb{Y}}_{p_0}'$ for which $P_s = P_w = (p_0, \infty]$.
	
	\subsection{} In contrast to the former case the spaces we now construct, apart from the basic property $P_s^c = P_w^c = [1, \infty]$, satisfy $P_s \varsubsetneq P_w$. Namely, for any fixed $p_0 \in [1, \infty)$ we construct a space $\widetilde{\mathbb{Y}}_{p_0}$ for which $P_s = (p_0, \infty]$ and $P_w = [p_0, \infty]$. We consider the cases $p_0 = 1$ and $p_0 > 1$ separately, similarly as it was done in Section 2.
	
	By $\widetilde{\mathbb{Y}}_1$ we denote the second generation space $(Y_\tau, \rho, \mu)$ with construction based on $\tau_n = n$ and $F(n,i) = 2^i$. Recall that $\mu$ is non-doubling.
	
	\begin{proposition}
		Let $\widetilde{\mathbb{Y}}_1$ be the metric measure space defined above. Then the associated centered operator $M^c$ is of strong type $(1,1)$, while the non-centered $M$ is of weak type $(1, 1)$, but is not of strong type $(1,1)$.
	\end{proposition}
	
	\begin{prof*}
		First note that it is easy to verify that $M^c$ is of strong type $(1,1)$, by using exactly the same argument as in the proof of Proposition 5.
		In the next step we show that $M$ is not of strong type $(1,1)$. Indeed, let $f_n = \delta_{y_n}$, $n \geq 1$. Then $\|f_n\|_1=d_n$ and for $i=1, \dots, n$ we have $Mf_n(y_{ni}') \geq (1+ 1 / n +2^i)^{-1} > 1 / 2^{i+1}$ and hence we obtain $\|Mf_n\|_1 \geq \sum_{i=1}^n 2^i d_n / 2^{i+1} = n \|f_n\|_1 / 2$.
		
		To complete the proof, it suffices to show that $M$ is of weak type $(1,1)$. Let $f \in \ell^1(\widetilde{\mathbb{Y}}_1)$, $f \geq 0$, and consider $\lambda > 0$ such that $E_\lambda (Mf) \neq \emptyset$. If $\lambda < A_{Y_\tau}(f)$, then $\lambda |E_\lambda(Mf)| / \|f\|_1 < 1$ follows. Therefore, from now on assume that $\lambda \geq A_{Y_\tau}(f)$.
		With this assumption, if for some $y \in T_n$ we have $Mf(y) > \lambda$, then any ball $B$ containing $y$ and realizing $A_B(f) > \lambda$ must be a subset of $T_n$. Take any $n \in \mathbb{N}$ such that $E_\lambda(Mf) \cap T_n \neq \emptyset$. If $\lambda < A_{T_n}(f)$, then
		\begin{equation} \label{3.2.1}
		\frac{\lambda |E_\lambda(Mf) \cap T_n|}{\sum_{y \in T_n} f(y)|\{y\}|} \leq 2.
		\end{equation}
		Assume that $\lambda \geq A_{T_n}(f)$ and take $y \in E_\lambda(Mf) \cap T_n$. Now, any ball $B$ containing $y$ and realizing $A_B(f) > \lambda$ must be a proper subset of $T_n$. First, consider the case $E_\lambda(Mf) \cap T_n' = \emptyset$. If $y_n \in E_\lambda(Mf) \cap T_n$, then we obtain $\sum_{y \in T_n \setminus T_n'} f(y)|\{y\}| > \lambda |\{y_n\}|$ and since $|E_\lambda(Mf) \cap T_n| \leq 2|\{y_n\}|$, (\ref{3.2.1}) follows. Otherwise, if $y_n \notin E_\lambda(Mf) \cap T_n$, then, 
		necessarily, $f(y) > \lambda$ for every $y \in E_\lambda(Mf) \cap T_n$ and hence (\ref{3.2.1}) again follows. Finally, in the case $E_\lambda(Mf) \cap T_n' \neq \emptyset$, denote $j = \max\{i \in \{1, \dots, n\} \colon Mf(y_{ni}')>\lambda\}$. Therefore, $\sum_{y \in T_n}f(y)|\{y\}| > \lambda |\{y_{nj}'\}|$ and combining this with the estimate $|E_\lambda(Mf) \cap T_n| \leq 2 |\{y_{nj}'\}|$, we conclude that (\ref{3.2.1}) follows. Since $\lambda |E_\lambda(Mf) \cap T_n| / \sum_{y \in T_n} f(y)|\{y\}| \leq 2$ for any branch $T_n$ such that $E_\lambda(Mf) \cap T_n \neq \emptyset$, we have $\lambda |E_\lambda(Mf)| / \|f\|_1 \leq 2$ and, consequently, $\|Mf\|_{1, \infty} \leq 2\|f\|_1$.
		$\raggedright \hfill \qed$
	\end{prof*} 
	
	Now, fix $p_0 \in (1, \infty)$ and consider $\widetilde{\mathbb{Y}}_{p_0} = (Y_\tau, \rho, \mu)$ with construction based on $\tau_n = \tau_n(p_0)$ and $F(n,i)= F_{p_0}(n, i)$, defined in the same way as in Section $2.2$, by using the auxiliary sequences $c_n$, $e_n$ and $m_{nj}$, $s_{nj}$, $j \in \{1, \dots, e_n\}$, $n \in \mathbb{N}$.

	\begin{proposition}
		Let $\widetilde{\mathbb{Y}}_{p_0}$ be the metric measure space defined above. Then the associated centered maximal operator $M^c$ is of strong type $(1,1)$, while the non-centered $M$ is of weak type $(p_0,p_0)$, but is not of strong type $(p_0,p_0)$.
	\end{proposition}
	
	\begin{prof*}
		First note once again that it is easy to verify that $M^c$ is of strong type $(1,1)$, by using the same argument as in the proof of Proposition 5. In the next step we show that $M$ is not of strong type $(p_0,p_0)$. Indeed, let $f_n = \delta_{x_n}$, $n \geq 1$. Then $\|f_n\|_{p_0}^{p_0}=d_n$ and for $i=1, \dots, \tau_n$ we have
		$Mf_n(y_{ni}') \geq (1+1/\tau_n+m_{nj(n,i)})^{-1} \geq (2(1+m_{nj(n,i)}))^{-1}$ and hence
		\begin{align*}
		\|Mf_n\|_{p_0}^{p_0} & \geq \sum_{j=1}^{e_n} \sum_{k=1}^{s_{nj}}\frac{d_n m_{nj}}{\big(2(1+m_{nj})\big)^{p_0}}  = d_n \sum_{j=1}^{e_n} \frac{s_{nj}m_{nj}}{\big(2(1+m_{nj})\big)^{p_0}} \\ & \geq d_n \sum_{j=1}^{e_n} \frac{2^{1-j-p_0} n c_n }{(1+m_{nj})^{p_0}} =  2^{-p_0} d_n \sum_{j=1}^{e_n} \frac{nc_n}{(1+n)^{p_0}} = 2^{-p_0} e_n \frac{nc_n}{(1+n)^{p_0}} \|f_n\|_{p_0}^{p_0}.
		\end{align*}
		Since $\lim_{n \rightarrow \infty} e_n = \infty$ and $\lim_{n \rightarrow \infty} \frac{n c_n}{(1+n)^{p_0}} = 1$, we are done.
		
		To complete the proof, it suffices to show that $M$ is of weak type $(p_0,p_0)$. Let $f \in \ell^{p_0}(\widetilde{\mathbb{Y}}_{p_0})$, $f \geq 0$, and consider $\lambda > 0$ such that $E_\lambda (Mf) \neq \emptyset$. If $\lambda < A_{Y_\tau}(f)$, then using the inequality $\|f\|_1 \leq \|f\|_{p_0} |Y_\tau|^{1/q_0}$, we obtain $ \lambda^{p_0} |E_\lambda(Mf)| / \|f\|_{p_0}^{p_0} < 1$. Therefore, from now on assume that $\lambda \geq A_{Y_\tau}(f)$. Take any $n \in \mathbb{N}$ such that $E_\lambda(Mf) \cap T_n \neq \emptyset$. If $\lambda < A_{T_n}(f)$, then using similar argument as above we have
		\begin{equation} \label{3.2.2}
		\frac{\lambda^{p_0} |E_\lambda(Mf) \cap T_n|}{\sum_{y \in T_n} f(y)^{p_0} |\{y\}|} \leq 1.
		\end{equation}
		Consider $\lambda \geq A_{T_n}(f)$. Assume that $E_\lambda(Mf) \cap T_n' = \emptyset$. If $\lambda < A_{T_n \setminus T_n'}(f)$, then (\ref{3.2.2}) again follows. Otherwise, if $\lambda \geq A_{T_n \setminus T_n'}(f)$, then we consider two cases. If $y_n \in E_\lambda(Mf)$, then we obtain $f(y_n) \geq \lambda$ and hence 
		\begin{displaymath}
		\frac{\lambda^{p_0} |E_\lambda(Mf) \cap T_n|}{\sum_{y \in T_n} f(y)^{p_0} |\{y\}|} \leq \frac{2 \lambda^{p_0} |\{y_n\}|}{\sum_{y \in T_n} f(y)^{p_0} |\{y\}|} \leq 2.
		\end{displaymath}
		In the other case, if $y_n \notin E_\lambda(Mf)$, then $f(y) > \lambda$ holds for every $y \in E_\lambda(Mf) \cap T_n$ and hence (\ref{3.2.2}) follows one more time. Now assume that $E_\lambda(Mf) \cap T_n' \neq \emptyset$. See that $|E_\lambda(Mf) \cap T_n| \leq 3|E_\lambda(Mf) \cap T_n'|$. Consider the case $f(y_n)< (1+1/\tau_n + m_{ne_n})\lambda / 3$. If $y_{ni}' \in E_\lambda(Mf) \cap T_n'$ for some $i \in \{1, \dots, \tau_n\}$, then  $f(y_{ni}') \geq \lambda / 3$ or $f(y_{ni})|\{y_{ni}\}| \geq |\{y_{ni}'\}| \lambda / 3$ and hence $f(y_{ni}')^{p_0}|\{y_{ni}'\}| + f(y_{ni})^{p_0}|\{y_{ni}\}| \geq  |\{y_{ni}'\}| (\lambda / 3)^{p_0}$, which implies
		\begin{displaymath}
		\frac{\lambda^{p_0} |E_\lambda(Mf) \cap T_n|}{\sum_{y \in T_n} f(y)^{p_0} |\{y\}|} \leq \frac{3 \lambda^{p_0} |E_\lambda(Mf) \cap T_n'|}{\sum_{y \in T_n} f(y)^{p_0} |\{y\}|} \leq 3^{p_0+1}.
		\end{displaymath}
		Finally, in the case $f(y_n) \geq (1+1/ \tau_n + m_{ne_n})\lambda / 3 $, denote $r = \min\{j \in \{1, \dots, e_n\} \colon f(y_n) \geq \frac{(1+1/ \tau_n+ m_{nj})\lambda}{3}\}$. Let $T_n^{(r)} = \{y_{ni}' \colon i \in \{1, \dots, \sum_{j=1}^{r-1}s_{nj}\}\}$. Note that the case $T_n^{(r)} = \emptyset$ is possible. Assume that $T_n^{(r)} \neq \emptyset$. If $y_{ni}' \in E_\lambda(Mf) \cap T_n^{(r)}$, then $f(y_{ni}')^{p_0}|\{y_{ni}'\}| + f(y_{ni})^{p_0}|\{y_{ni}\}| \geq  |\{y_{ni}'\}| (\lambda / 3)^{p_0} $ and hence
		\begin{displaymath}
		\frac{\lambda^{p_0} |E_\lambda(Mf) \cap T_n^{(r)}|}{\sum_{i \colon y_{ni}' \in T_n^{(r)}} (f(y_{ni}')^{p_0}|\{y_{ni}'\}| + f(y_{ni})^{p_0}|\{y_{ni}\}|)} \leq 3^{p_0+1}.
		\end{displaymath}
		Moreover, we have
		\begin{align*}
		\frac{\lambda^{p_0} |E_\lambda(Mf) \cap \big(T_n \setminus T_n^{(r)}\big)|}{f(y_n)^{p_0} |\{y_n\}|} & \leq \Big(\frac{3}{1+ m_{nr}}\Big)^{p_0} \frac{|T_n \setminus T_n^{(r)}|}{|\{y_n\}|} \leq \Big(\frac{3}{1+m_{nr}}\Big)^{p_0} \frac{3 |(T_n \setminus T_n^{(r)}) \cap T_n')|}{|\{y_n\}|} \\
		& \leq \Big(\frac{3}{1+m_{nr}}\Big)^{p_0} 3 \sum_{j=r}^{e_n}n c_n 2^{2-j} < 2^{3-r}3^{p_0+1} n c_n \Big(\frac{1}{1+m_{nr}}\Big)^{p_0} \\ & = 4 \cdot 3^{p_0+1} \frac{nc_n}{(1+n)^{p_0}} \leq 4 \cdot 3^{p_0+1}.
		\end{align*}
		Therefore, regardless of the possibilities, $T_n^{(r)} = \emptyset$ or $T_n^{(r)} \neq \emptyset$, we obtain $\lambda^{p_0} |E_\lambda(Mf) \cap T_n| / \sum_{y \in T_n} f(y)^{p_0} |\{y\}| \leq 4 \cdot 3^{p_0+1}$. Since $ \lambda |E_\lambda(Mf) \cap T_n| / \sum_{y \in T_n} f(y)^{p_0} |\{y\}| \leq 4 \cdot 3^{p_0+1}$ for any branch $T_n$ such that $E_\lambda(Mf) \cap T_n \neq \emptyset$, we have $\lambda^{p_0} |E_\lambda(Mf)| / \|f\|_{p_0}^{p_0} \leq 4 \cdot 3^{p_0+1}$ and consequently $\|Mf\|_{p_0, \infty}^{p_0} \leq 4 \cdot 3^{p_0+1} \|f\|_{p_0}^{p_0}$. 
		$\raggedright \hfill \qed$
	\end{prof*}
	
	\section{Proof of Theorem 1}
	
	All spaces discussed above were constructed in such a way as to avoid any interactions between the different branches in the context of considerations relating to the existence of the weak and strong type inequalities. Therefore we can construct a new space consisting of two types of branches, one borrowed from some first generation space and one from some second generation space, and to ensure that the operators $M^c$ and $M$ inherit a particular property of a particular space. We explain the construction of such a space in detail proving Theorem 1.
	
	\begin{proof*}
		We consider a few cases. If the equalities $P_s^c = P_s$ and $P_w^c = P_w$ are supposed to hold, then the expected space may be chosen to be a first generation space. If, in turn, we have $P_s^c = P_w^c = [1, \infty]$, but $P_s \neq [1, \infty]$, then the expected space may be chosen to be a second generation space. Finally, in other cases we can find spaces $\mathbb{X} = (X, \rho_X, \mu_X)$ and $\mathbb{Y} = (Y, \rho_Y, \mu_Y)$, of first and second generation, respectively, for which
		\begin{itemize}
			\item $P_s^c(\mathbb{X}) = P_s(\mathbb{X}) = P_s^c$ and $P_w^c(\mathbb{X}) = P_w(\mathbb{X}) = P_w^c$, \smallskip
			\item $P_s^c(\mathbb{Y}) = P_w^c(\mathbb{Y}) = [1, \infty]$, $P_s(\mathbb{Y}) = P_s$ and $P_w(\mathbb{Y}) = P_w$. \smallskip
		\end{itemize}
		Using $\mathbb{X}$ and $\mathbb{Y}$ and assuming that $X \cap Y = \emptyset$ we construct the space $\mathbb{Z} = (Z, \rho_Z, \mu_Z)$ as follows.
		Denote $Z = X \cup Y$. We define the metric $\rho_Z$ on $Z$ by
		\begin{displaymath}
		\rho_Z(x,y) = \left\{ \begin{array}{rl}
		\rho_X(x,y) & \textrm{if }  \{x,y\} \subset X,   \\
		\rho_Y(x,y) & \textrm{if }  \{x,y\} \subset Y,   \\
		2 & \textrm{in the other case,} \end{array} \right. 
		\end{displaymath} 
		and the measure $\mu_Z$ on $Z$ by
		\begin{displaymath}
		\mu_Z(E) = \mu_X(E \cap X) + \mu_Y(E \cap Y), \qquad E \subset Z.
		\end{displaymath}
		It is not hard to show that $\mathbb{Z}$ has the following properties
		\begin{itemize}
			\item $P_s^c(\mathbb{Z})=P_s^c(\mathbb{X}) \cap P_s^c(\mathbb{Y})= P_s^c \cap [1, \infty] = P_s^c$, \smallskip
			\item $P_s(\mathbb{Z}) = P_s(\mathbb{X}) \cap P_s(\mathbb{Y}) = P_s^c \cap P_s = P_s$, \smallskip
			\item $P_w^c(\mathbb{Z})=P_w^c(\mathbb{X}) \cap P_w^c(\mathbb{Y})= P_w^c \cap [1, \infty] = P_w^c$, \smallskip
			\item $P_w(\mathbb{Z}) = P_w(\mathbb{X}) \cap P_w(\mathbb{Y}) = P_w^c \cap P_w = P_w$, 
		\end{itemize} 
		and therefore it may be chosen to be the expected space.
		Finally, it is not hard to see that $\mu_Z$ is non-doubling, since it is bounded and there is a ball $B$ in $Z$ with radius $r = 1$ and $|B| < \epsilon$ for any arbitrarily small $\epsilon > 0$.
		$\raggedright \hfill \qed$
	\end{proof*}
	
	\section*{Acknowledgement}
	This article was largely inspired by the suggestions of my supervisor Professor Krzysztof Stempak. I would like to thank him for insightful comments and continuous help during the preparation of the paper.

\end{document}